\documentclass[12pt,twoside]{amsart}
\usepackage{amsmath}
\usepackage{epsfig}
\title{A new construction of wavelet sets} 
\author{ Eugen J.\ Ionascu}
\curraddr{Columbus State University\\4225 University Avenue\\ Columbus, GA 31907}
\email{ionascu\_eugen@\ colstate.edu}
\subjclass{Primary 42C40, 42C30; Secondary 03E15}
\date{August 30th, 2002}
\keywords{wavelet set, wavelet induced isomorphism, Schr\"{o}der-Cantor-Bernstein construction} 
\begin{document}
\baselineskip18pt 
                               
\def\RR{{\rm I}\!{\rm R}}
\def\fp#1{(#1)}
\newtheorem{theorem}{\hspace{\parindent}
T{\scriptsize HEOREM}}[section]
\newtheorem{proposition}[theorem]
{\hspace{\parindent }P{\scriptsize ROPOSITION}}
\newtheorem{corollary}[theorem]
{\hspace{\parindent }C{\scriptsize OROLLARY}}
\newtheorem{lemma}[theorem]
{\hspace{\parindent }L{\scriptsize EMMA}}
\newtheorem{definition}[theorem]
{\hspace{\parindent }D{\scriptsize EFINITION}}
\newtheorem{problem}[theorem]
{\hspace{\parindent }P{\scriptsize ROBLEM}}
\newtheorem{conjecture}[theorem]
{\hspace{\parindent }C{\scriptsize ONJECTURE}}
\newtheorem{example}[theorem]
{\hspace{\parindent }E{\scriptsize XAMPLE}}
\newtheorem{remark}[theorem]
{\hspace{\parindent }R{\scriptsize EMARK}}
\renewcommand{\thetheorem}{\arabic{section}.\arabic{theorem}}
\newcommand{\lr}{L^2(\RR)} 
\newcommand{\du}{\stackrel{.}{\bigcup}}
\renewcommand{\theenumi}{(\roman{enumi})}
\renewcommand{\labelenumi}{\theenumi} 
\def\RR{{\rm I}\!{\rm R}}
\def\LL{{\rm I}\!{\rm L}}
\def\NN{{\rm I}\!{\rm N}}
\def\MM{{\rm I}\!{\rm M}}
\def\QQ{{\rm I}\!\!\!{\rm Q}}
\def\ZZ{{\rm Z}\!\!{\rm Z}}
\def\CC{{\rm I}\!\!\!{\rm C}}
\def\vp{\varepsilon}
\def\phi{\varphi}
\def\ra{\rightarrow}
\def\sd{\bigtriangledown}
\def\ac{\mathaccent94}
\def\wi{\sim}
\def\wt{\widetilde} 
\def\bb#1{{\Bbb#1}}
\def\bs{\backslash}
\def\cal{\mathcal}
\def\ca#1{{\cal#1}}
\def\Bbb#1{\bf#1}
\def\blacksquare{{\ \vrule height7pt width7pt depth0pt}} 
\def\bsq{\blacksquare}
\def\proof{\hspace{\parindent}{P{\scriptsize ROOF}}}
\def\pofthe{P{\scriptsize ROOF OF} 
T{\scriptsize HEOREM}\  }
\def\pofle{\hspace{\parindent}P{\scriptsize ROOF OF} 
L{\scriptsize EMMA}\  }
\def\pofcor{\hspace{\parindent}P{\scriptsize ROOF OF} 
C{\scriptsize ROLLARY}\  }
\def\pofpro{\hspace{\parindent}P{\scriptsize ROOF OF} 
P{\scriptsize ROPOSITION}\  }
\def\n{\noindent}
\def\iter#1{^{(#1)}}
\def\bh{\ca B(H)}
\def\ld{\overline}
\def\vpe{V_{\psi}^{\eta}}
\def\wu{\ca W(\ca U)}
\def \dper {\ca D\ca P}
\def\lir{L^{\infty}(\RR)} \def\ltr{L^2(\RR)} \def\ltrn{L^2(\RR^n)}
\def\lie{L^{\infty}(E)}\def\lte{L^2(E)}
\def\mpe{\ca M_{\psi,\eta}}
\def\wh{\widehat}
\def\eproof{$\hfill\bsq$\par}
\def\ws{\ca W\ca S}
\def\wi{\ca W\ca I}
\def\dst{\displaystyle}
\def\wsna{{\ca W}{\ca S}(n,A) }
\def\daw{{\ca W}_A}
\def\fwi{\ca W\ca I_1}
\def\swi{\ca W\ca I_2}
\def\ds{\displaystyle}
\begin{abstract}                                                        
We show that the class of (dyadic) wavelet sets is in one-to-one correspondence 
to a special class of Lebesgue measurable isomorphisms of $[0,1)$ which we call {\it wavelet induced} maps.  We then define two natural classes of maps 
$\fwi$ and $\swi$ which, in order to simplify their construction, retain only part of the characterization properties of a wavelet induced map. We prove that 
each wavelet induced map appears from the Schr\"{o}der-Cantor-Bernstein construction applied to some $u\in \fwi$  and $v\in \swi$. Consequently, the 
construction of a wavelet set is basically equivalent to the easier construction of two maps $u\in \fwi$  and $v\in \swi$. Some older results on wavelet sets are recovered using this new point of view.  The connectivity result of Speegle (\cite{sp}) is recaptured and the completeness in the natural metric of the class of wavelet sets is reestablished. Although these ideas seem to generalize to more than one dimension, specific examples are given only in the one dimensional case.
\end{abstract}                                                    
\maketitle                                                       
\section{Introduction}
In \cite{dl} the authors introduced the notion of wavelet set which turned out to be one of the building blocks of their approach to wavelet analysis from an operator theory point of view. At the same time, and independently, the notion of wavelet set appeared as the support set of so called MSF-wavelets (minimally supported frequency) in a series of papers: \cite{fw}, \cite{hww1} and \cite{hww2}. 

One easy way to fabricate a wavelet is to normalize the Fourier transform of the characteristic function of a wavelet set. Wavelet sets have been generalized to $n$ dimensions (see \cite{dls1} and \cite{dls2}). The important result of the existence of wavelets for unitary systems having an expansive dilation matrix was based on the existence of wavelet sets. These ideas were taken into the realm of frame theory and the notion of wavelet set was generalized even further to frame (tight frame or normalized tight frame)  wavelet sets in \cite{ddg1}, \cite{ddg2}, \cite{ddgh1}-\cite{ddgh3}. In \cite{bmm} the authors give an ingenious description of how one can construct a wavelet set. The purpose of this paper is to consider a different approach to the construction of (dyadic) wavelet sets which is purely set theoretic. 

\section{Preliminary Results}
We denote by $\mu$  the Lebesgue measure on $\RR$. 
The $L^2$-space with respect to $\mu $ will be written simply as  $L^2(\RR)$.  An {\it orthonormal} wavelet is (cf. \cite{dl}) a function $w\in L^2(\RR)$ for which the family of functions $\{w_{j,k}\}_{j,k\in \ZZ}$  defined by 
\begin{equation}
w_{j,k}(s)= 2^{j/2}w(2^js-k),\ \ s\in \RR, \ j, k\in \ZZ,
\end{equation}
is an orthonormal basis for $\lr$.\par
We say that a measurable subset $W$ of $ \RR$ is a {\it wavelet set}  if $\frac{1}{\sqrt{\mu(W)}}\chi_W=\wh{w}$,
where $w$ is a wavelet in $L^2(\RR)$ 
and $\wh{w}$ is  the Fourier-Plancherel transform on $L^2(\RR)$ of the function $w$ and which for $f\in L^1(\RR)\cap \lr$ is defined by                                   
$$\wh f (x)=\frac{1}{\sqrt{2\pi}}\int_{\RR}                            
e^{-itx}f(t)d\mu (t),\ \  \ x\in \RR.$$ 
\par

One of the simplest examples of wavelet sets is the 
Littlewood-Paley wavelet set $E:=[-2\pi,-\pi)\cup[\pi,2\pi)$.
A less obvious example is the following union of eight intervals 
\begin{equation}\label{eq:S}
\dst S:=\left.                           
\begin{array}{l}
\dst \left[-\frac{4\pi}{3},- \frac{5\pi}{4}\right)\cup
\left[-\pi,- \frac{2\pi}{3}\right)\cup
\left[-\frac{5\pi}{8},- \frac{\pi}{2}\right)\cup 
\left[\frac{4\pi}{7}, \frac{2 \pi}{3}\right)\\ \\
\dst \cup \left[\frac{3\pi}{4},\pi\right)\cup 
\left[\frac{4\pi}{3},\frac{11\pi}{8}\right)
\cup \left[4\pi, \frac{32 \pi}{7}\right)
\cup \left[\frac{11\pi}{2}, 6\pi \right).
\end{array}  
\right.
\end{equation}

The next result was announced independently in 
\cite{fw} and \cite{dl} and it is definitely 
the first step in a better understanding of the notion of a wavelet set. 
We refer the reader  to \cite{ilp} for a proof of this proposition.
In order to state the result let us 
introduce some notation. Let $\tau :\RR\ra E$ be the function                   
defined by $\tau(x)=x+2j\pi$, where $j$ is             
the unique integer satisfying $x +2j\pi\in E$         
and let $\delta:\RR\bs \{0\}\ra E$ be the map     
defined by $\delta(x)=2^kx$, where $k$ is the unique   
integer for which $2^k x\in E$.

\begin{proposition} The following conditions are equivalent  for any measurable subset $W$ of $\RR$ $:$\par
\begin{enumerate} 
\item $W$ is a wavelet  set, \par 
\item there exists a set $W'$ such that $W'=W$ a.e., the family of sets $\{W'+2k\pi\}_{k\in\ZZ}$ is a partition of $\RR$ and, at the same time, the family $\{2^kW'\}_{k\in \ZZ}$ is a partition of $\RR\setminus\{0\}$,
\item  there exists a set $W''$ such that $W''=W$ a.e. and  $\tau_{|W''},\delta_{|W''}:W''\ra E$  are measurable bijections.
\end{enumerate}
\end{proposition}

In \cite{ilp} a wavelet set having the property of $W'$ in (ii) (or equivalently the property of $W''$ in (iii))
was called {\it regularized}. Let us denote by $\ws$  the class of all 
wavelet sets.
 
The class $\ws$ is very rich. In \cite{ip} it was shown that every point
$x_0\in \RR\setminus\{0\}$ contains a neighborhood which is a part of a wavelet set. In \cite{sp} it was proved that $\ws$ is path-connected (in the norm topology on $L^2(\RR)$ when $\ws$ is naturally imbedded in $L^2(\RR)$). It was shown in 
\cite{gs} that $\ws$ becomes a complete metric space $(\ws,d)$ with the  metric:
\begin{equation}\label{distance}
d(W_1,W_2):=\mu(W_1\bigtriangledown W_2)^{\frac{1}{2}} +
\left(\int_ {W_2\bigtriangledown W_1} \frac{1}{|x|} d\mu(x)\right)^{\frac{1}{2}},
\end{equation}
where $W_1\bigtriangledown W_2=(W_1\setminus W_2)\cup (W_2\setminus W_1)$.
In spite of this richness, it is not very obvious how would one construct a wavelet set. 

Thus by part (iii) of this proposition, we can associate with every $W\in \ws$ a measurable bijection on $E$  defined by 
\begin{equation}\label{eq:hf}
h_{W}:= \tau_{| W'}\circ (\delta _{| W'})^{-1}.
\end{equation}
By its definition, this map is essentially uniquely determined by 
$W$ in the sense that for two sets  $W'$ 
and $W''$  as in Proposition~2.1, the maps $h_{W'}$, $h_{W''}$ 
coincide almost everywhere [$\mu$]. 
Thus we will denote this map simply by $h_W$.
It turns out that the conjugation 
$\wt h_W:=\xi \circ h_W\circ \xi ^{-1}: [0,1)\ra [0,1)$
of $h_W$, by the function $\xi:E\ra [0,1)$ defined by                 
\begin{equation}\label{eq:xi}                  
\dst \xi(x)=\left \{                           
\begin{array}{l}                               
\dst \frac{x}{2\pi},\ \ x\in [\pi,2\pi),\\  \\  
\dst \frac{x}{2\pi}+1,\ \ x\in [-2\pi,-\pi ),  
\end{array}                                    
\right.                                        
\end{equation}
takes a simpler form than $h_W$. 
\begin{definition}
For every wavelet set $W$ the map $\wt h_W$ constructed as above is called an {\it wavelet induced} isomorphisms
of $[0,1)$.
\end{definition}
We have the following
characterization of the class of wavelet sets in terms 
of the corresponding maps $\wt h_W$.
\begin{proposition} \label{prop:charh}                               %
Let $W\in \ws$ and $\wt h_W$ be defined as above. Then               %
the map $\wt h_W$ has the following properties:                      %
\begin{enumerate}                                                    %
\item $\wt h_W$ is a measurable bijection of [0,1),                  %
\item  there exists a measurable partition                           %
$\{A_k\}_{k\in \ZZ}$ of  $[\frac{1}{2},1)$                           %
and a measurable partition                                           %
$\{B_k\}_{k\in \ZZ}$ of $[0,\frac{1}{2})$, such that                 %
\begin{equation}\label{eq:wthf}                                      %
\dst \wt h_W(x)=\left \{                                             %
\begin{array}{l}                                                     %
\dst \lfloor 2^kx \rfloor ,\ \ x\in A_k,\  k\in \ZZ,\\  \\           %
\dst \lfloor 2^k(x-1) \rfloor ,\ \ x\in B_k,\ k\in \ZZ,              %
\end{array}                                                          %
\right.                                                              %
\end{equation}                                                       %
where $\lfloor x \rfloor$ denotes the fractional                     %
part of the real number $x$.                                         %
\item    if $h$ is a map satisfying $(i)$                            %
and $(ii)$ then there exists a wavelet set $W$ such $h=\wt h_W$;     
\item for two wavelet sets $W_1$ and $W_2$ we have 
$W_1\cap W_2=\delta_{|W_1}^{-1}\circ \xi^{-1}
(\Omega)=
\delta_{|W_2}^{-1}\circ \xi^{-1}(\Omega)$, a.e. $[\mu]$, where  
$\Omega=\{x\in [0,1): \wt h_{W_1}(x)=\wt h_{W_2}(x)\}$.       %
\end{enumerate}                                                      %
\end{proposition}                                                    %
\proof.\  One can easily check that $\tau(t)= \xi^{-1}(\lfloor \frac{t}{2\pi} \rfloor)$ for every $t\in \RR$. First, let us observe that if $x\in [\frac{1}{2},1)$ then 
$u=\xi^{-1}(x)=2\pi x $ and let us write $\delta ^{-1}(u)=2^ku$ 
with $k\in \ZZ$. Hence, using the formula mentioned above, we have $\tau(2^ku)=\xi^{-1}(\lfloor \frac{2^ku}{2\pi} \rfloor)=\xi^{-1}(\lfloor 2^kx \rfloor)$.
Thus $\wt h_W(x)=\xi(\tau_{|W}( 2^ku))=\lfloor 2^kx \rfloor$.\par
If $x\in [0,1/2)$, then $u=\xi^{-1}(x)=2\pi (x-1)$
and if we write $\delta ^{-1}(u)=2^ku$ with $k\in \ZZ$, we get 
$\wt h_W(x)=\xi(\tau_{|W}( 2^ku))=\lfloor \frac{2^ku}{2\pi} \rfloor=
\lfloor 2^k(x-1) \rfloor$. This proves claim (ii) of the proposition.\par
To prove claim (iii), let us consider that $h$ has the properties $(i)$ and $(ii)$
and we denote $h_1=\xi^{-1}\circ h\circ \xi:E\ra E$. Because of (\ref{eq:wthf}),
we obtain that for every $x\in E$ there exist $k(x),\ l(x)\in \ZZ$ such that 
$h_1(x)=2^{k(x)} x+2l(x)\pi$. By the assumptions on $h$ the maps $x\rightarrow k(x)$,
$x\rightarrow l(x)$ are measurable. Then we define $\psi: E\ra \RR$ by $\psi(x)=2^{k(x)}x$, $x\in E$, and $W:=\psi(E)$.  Clearly, $W$ is a measurable set, $\psi$ is one-to-one and $(\delta_{|W})^{-1}=\psi$.
Finally we define $\phi: W\ra E$ by $\phi(y)=y+2l(\psi^{-1}(y))\pi$. Then one can check that $h_1=\phi\circ \psi$ and $\phi(y)=\tau_{|W}(y)$, $y\in W$.  Since $h_1$ is one-to-one and $\psi$ is onto we conclude that $\phi$ is one-to-one. 
Also, $\phi$ is onto since $h_1$ is.
By Proposition~2.1 we see that $W$ is a wavelet set. According to (\ref{eq:hf}) and using the above facts it follows that $h_W=h_1$ and so  $\wt h_W=h$.\par

To prove (iv) let us show the equality in question by double inclusion.
For $y\in W_1\cap W_2$ we have $s:=\delta_{|W_1}(y)=\delta_{|W_2}(y)=\delta(y)$. Then $y=\delta_{|W_1}^{-1}(s)=\delta_{|W_2}^{-1}(s)$. Let us denote $\xi(s)$ by $u$. In other words $s=\xi^{-1}(u)$. This lets us write $y=\delta_{|W_1}^{-1}\circ
\xi^{-1}(u)=\delta_{|W_2}^{-1}\circ \xi^{-1}(u)$. This means that we need to check that $u\in \{x\in [0,1): \wt h_{W_1}(x)=\wt h_{W_2}(x)\}$ or equivalently 
$u\in \{\xi(t): t\in E, \ \   h_{W_1}(t)= h_{W_2}(t)\}$. Since $u=\xi(s)$ we need to see why is it true that $h_{W_1}(s)= h_{W_2}(s)$. Using (\ref{eq:hf}) this last equality 
is the same as $\tau_{|W_1}(y)=\tau_{|W_2}(y)$ which is true. This argument shows that 
$W_1\cap W_2\subset \delta_{|W_1}^{-1}\circ \xi^{-1}(\Omega)$ and 
$W_1\cap W_2\subset \delta_{|W_2}^{-1}\circ \xi^{-1}(\Omega)$.
\par

For the opposite inclusion let us start with $y=\delta_{|W_1}^{-1}\circ \xi^{-1}(u)$
where $u$ satisfies $\wt h_{W_1}(u)=\wt h_{W_2}(u)$. As before we let $s:=\xi^{-1}(u)$. This implies $ h_{W_1}(s)= h_{W_2}(s)$ or $\tau_{|W_1}\circ \delta_{|W_1}^{-1}(s)=\tau_{|W_2}\circ \delta_{|W_2}^{-1}(s)$. Taking in account that 
$\tau(t_1)=\tau(t_2)$ implies $t_1=t_2+2k\pi$ for some $k\in \Bbb Z$ we obtain that $\delta_{|W_1}^{-1}(s)= \delta_{|W_2}^{-1}(s)+2k\pi$, $k\in \Bbb Z$. Let $\delta_{|W_1}^{-1}(s)=2^ns$ and $\delta_{|W_2}^{-1}(s)=2^ms$, for some $m$ and $n\in \Bbb Z$. If $n\not=m$ then $s=\frac{2k\pi}{2^n-2^m}$.
It is clear that the set 
\begin{equation}\label{calf}
\ca F=\left\{\ds \frac{2l\pi}{2^i-2^j}:i,j,l\in \Bbb Z,i\not=j\right\}
\end{equation}
 is countable and so it has measure zero. Assuming that $s\not \in \ca F$, we get $n=m$ and $k=0$ which implies $y=\delta_{|W_1}^{-1}\circ \xi^{-1}(u)=\delta_{|W_1}^{-1}(s)=\delta_{|W_2}^{-1}(s)\in W_1\cap W_2$. This argument proves that 
$\delta_{|W_1}^{-1}\circ \xi^{-1}
(\Omega)\setminus \ca F\subset W_1\cap W_2$ and 
$\delta_{|W_2}^{-1}\circ \xi^{-1}(\Omega)\setminus 
\ca F\subset W_1\cap W_2$.\eproof

For the  wavelet set $S$  defined by (2), we computed the 
map $\wt h_S$ and obtained
\begin{equation}\label{eq:example}
\dst \wt h_S(x)=\left \{       
\begin{array}{l}               
\dst \lfloor 2^{-1}(x-1) \rfloor \quad \text{on} \quad 
\left[0,\frac{1}{3}\right)\cup 
\left[\frac{3}{8},\frac{1}{2}\right),\\ \\
\dst \lfloor 2^2x \rfloor  \quad \text{on} \quad 
\left[\frac{1}{2}, \frac{4}{7}\right)\cup \left[\frac{11}{16},\frac{3}{4}\right),\\ \\
\dst \lfloor 2^{-1}x \rfloor  \quad \text{on}  \quad              
\left[\frac{4}{7}, \frac{2}{3}\right)\cup 
\left[\frac{3}{4},1\right),\\ \\  
\dst \lfloor x \rfloor  \quad \text{on} \quad 
\left[\frac{1}{3}, \frac{3}{8}\right)\cup 
\left[\frac{2}{3}, \frac{11}{16}\right).
\end{array} \right.                                          
\end{equation}  

One can check that this map is a measurable bijection from
$[0,1)$ into $[0,1)$.  The graph of it is on page 6. \par
\begin{center}
\epsfig{file={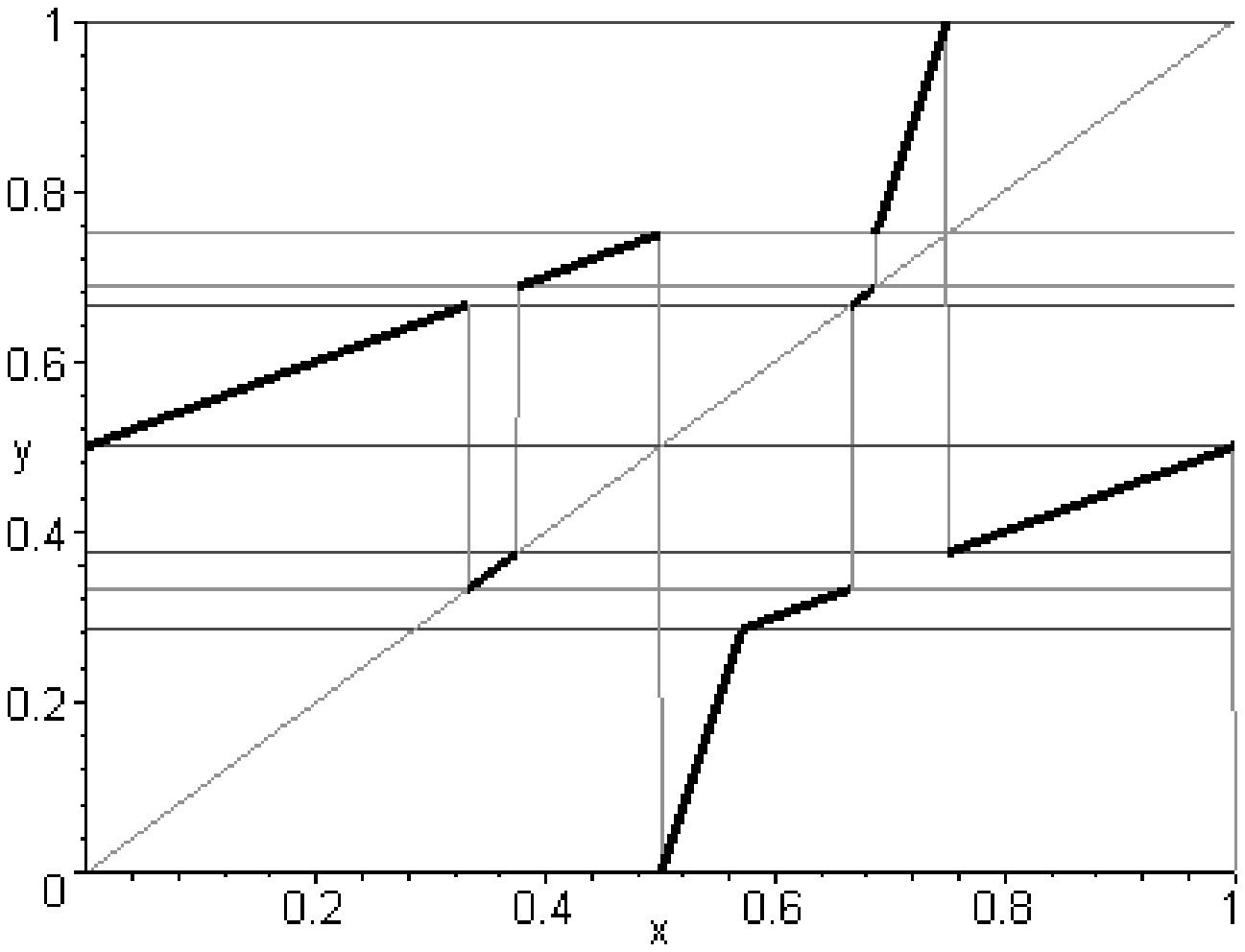},height={4in},width={4in}}
\end{center}

We denote the class of all wavelet induced isomorphisms by $\wi$. The wavelet induced map for
and wavelet set played an important role in the series of papers \cite{ai}, \cite{ailp}
and \cite{ilp}. One essential hypothesis that we needed in these works was the existence 
of a measurable cross section for the isomorphism (i.e. a measurable set which contains exactly one point from 
each orbit of the isomorphism). This existence was partially solved in \cite{ai} but it is still an open conjecture in the general situation. In this paper we have no need for such an assumption. However, we obtain a similar path of wavelets in Theorem~\ref{connect} as in \cite{ilp}.  By Proposition~\ref {prop:charh}  every  map in $\wi$ can be expressed as in (\ref {eq:wthf}).
\section{The Schr\"{o}der-Cantor-Bernstein construction}
Let $\fwi$ be the class of all measurable maps $f:[0,1)\ra [0,1)$ that are defined by a measurable partition
$\{A_k\}_{k\in \NN}$ of $[\frac{1}{2},1)$                            
and a measurable partition                                           
$\{B_k\}_{k\in \NN}$ of $[0,\frac{1}{2})$, such that                 
\begin{equation}\label{eq:cwi}                                      
\dst f(x)=\left \{                                             
\begin{array}{l}                                                   
\dst  \frac{x}{2^k} ,\ \ x\in A_k,\  k\ge 1,\\  \\            
\dst \frac {x-1}{2^k}+1,\ \ x\in B_k,\ k\ge 1.             
\end{array}                                                          
\right.                                                              
\end{equation}                                                       
The following is a simple consequence of the above definition.
\begin{lemma}\label{lem:firstm}
Every function $f\in \fwi$ is one-to-one and $\mu(f(\sigma))\le \frac{1}{2}\mu(\sigma)$ for every measurable subset $\sigma$ of $[0,1)$.
\end{lemma} 
\proof. \ Let $x_1,x_2\in [0,1)$, $x_1\ne x_2$. We need to analyze essentially four cases. If $x_1,\ x_2\in A_k$ or 
$x_1,\ x_2\in B_k$ for some $k\ge 1$ clearly $f(x_1)\ne f(x_2)$. If $x_1\in A_k$ and $x_2\in B_l$ for some $k,l\ge 1$ then $f(x_2)=(x_2-1)/2^l+1\ge 1/2>f(x_1)$.
Suppose now $x_1\in A_k$ and $x_2\in A_l$ for $k>l\ge 1$. In this case $x_1,\ x_2\ge 1/2$ and so $f(x_2)=x_2/2^l\ge 1/2^{l+1}\ge 1/2^k>f(x_1)$. Finally, if $x_1\in B_k$ and $x_2\in B_l$ for some $k>l\ge 1$ we have $x_1,\ x_2\in [0,1/2)$. Hence, $f(x_2)=(x_2-1)/2^l+1<
1-1/2^{l+1}\le 1-1/2^k\le f(x_1)$.\par
For the second part of this lemma let us observe that we can write 
$\sigma$ as a disjoint union $\sigma=\bigcup_{n\ge 1} \sigma_n \cup \bigcup_{n\ge 1}\sigma_n'$ where $\sigma_n=\sigma\cap A_n$, $\sigma_n'=\sigma\cap B_n$, $n\in \NN$. 
Clearly, $\mu(f(\sigma_n))\le (1/2)\mu(\sigma_n)$ and  
$\mu(f(\sigma_n'))\le (1/2)\mu(\sigma_n')$ for all $n\in \NN$. Adding up all these 
inequalities we obtain that $\mu(f(\sigma))\le \mu(\sigma)$ for every measurable set 
$\sigma$.\eproof

{\bf Remark.} Let us emphasize the fact that constructing a map $f$ in $\fwi$ is just a simple matter of choosing two measurable partitions: one for $[0,1/2)$ and one for $[1/2,1)$. The next figure shows the line segments of equations involved in (\ref{eq:cwi}) and gives an idea of why Lemma~\ref{lem:firstm} is true. 
\begin{center}
\epsfig{file={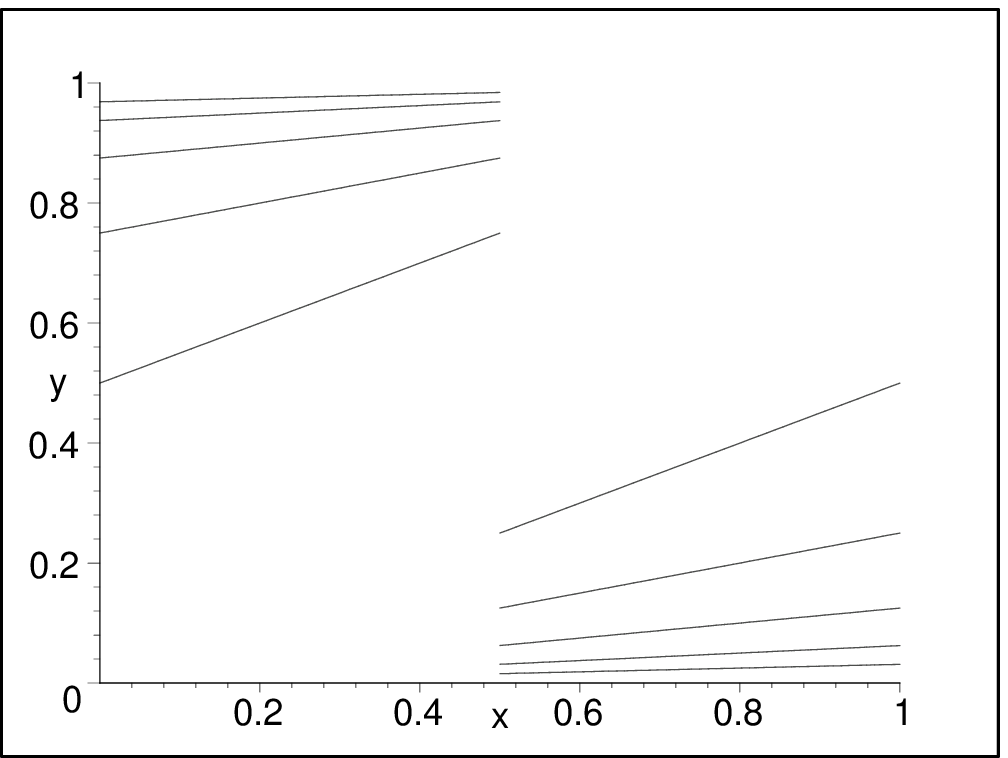},height={4in},width={4in}}
\end{center}

We introduce now the class $\swi$ of all measurable one-to-one maps 
$g:[0,1)\ra [0,1)$ so that for every $x\in [0,1)$ there exist $k,l\in \ZZ, k,l\ge 0$ such that 
$\dst g(x)=\frac{x+l}{2^k}$. 
We remind the reader the following fact from set theory known as the 
Schr\"{o}der-Cantor-Bernstein theorem.
\begin{proposition}$(${\bf Schr\"{o}der-Cantor-Bernstein}$)$ 
\label{prop:cabe}
Let $A$ and $B$ two arbitrary sets, $u:A\ra B$ and $v:B\ra A$ be two one-to-one maps. 
Then the map $u\diamond v:A\ra B$ defined by  
\begin{equation}\label{cabe}
\dst (u\diamond v )(x)=\left \{       
\begin{array}{l}               
\dst u(x) \quad \text{for} \quad x\in 
\bigcup_{k=0}^{\infty}(v\circ u)^k(A\setminus v(B)),\\ \\  
\dst v^{-1}(x)  \quad \text{for} \quad x\in 
\bigcup_{k=0}^{\infty}(v\circ u)^k(v(B)\setminus (v\circ u) (A))\cup 
\bigcap_{k=0}^{\infty}(v\circ u)^k(A),
\end{array} \right.                                          
\end{equation}  
is a bijection.
\end{proposition}
\n {\bf Remark.} It is easy to see that the inverse function of $u\diamond v$ given as in 
(\ref{cabe}) is in fact $v\diamond u$.
We are now ready for the main result of this note. This next theorem has almost the flavor of a factorization theorem. 
\begin{theorem}\label{main}
Every wavelet induced isomorphism $h$ is the result of the Schr\"{o}der-Cantor-Bernstein 
construction, i.e. $h=u\diamond v$ a.e. for some $u\in \fwi$ and $v\in \swi$ where $\diamond$ is defined in (\ref{cabe}).
Conversely, every map $u\diamond v$ with $u\in \fwi$ and $v\in \swi$ is an wavelet induced isomorphism. (The writing  $h=u\diamond v$ is in general not unique.)
\end{theorem} 
\proof.\ The last part of the theorem follows from the Proposition~2.2 and the fact that 
$u\diamond v$ is by construction a measurable bijection of the form (ii) in Proposition~2.2.
To show the first part, let us start with $h\in \wi$ and define $\ca D_1$ be the measurable set of all $x\in [0,1)$ for which the $k$ in the definition (\ref{eq:wthf})  is a negative integer. We let $u(x)=h(x)$ for $x\in \ca D_1$ and extend $u$ to $\ca D_2=[0,1)\setminus \ca D_1$ such that $ u\in \fwi$. This extension is not unique but it can be easily constructed as observed in the remark following the Lemma~\ref{lem:firstm}.
(For instance, for the map given in (\ref{eq:example})
we can take
\begin{equation}\label{eq:f}
\dst u(x)=\left \{       
\begin{array}{l}               
\dst \frac{x-1}{2}+1 \quad \text{on} \quad 
\left[0,\frac{1}{2}\right),\\ \\
\dst \frac{x}{2}  \quad \text{on} \quad 
\left[\frac{1}{2}, 1\right),
\end{array} \right.                                          
\end{equation}  
in which case $\ca D_1=\left[0,\frac{1}{3}\right)\cup 
\left[\frac{3}{8},\frac{1}{2}\right)\cup \left[\frac{4}{7}, \frac{2}{3}\right)\cup 
\left[\frac{3}{4},1\right)$.) 
To continue the proof, let $\ca R_1:=u(\ca D_1)$ and $\ca R_2:=[0,1)\setminus \ca R_1$. Next we define $v$ on $\ca R_2$ as $v(x)=h^{-1}(x)$ ($x\in \ca R_2$). Let us observe that  $v$ is defined in accordance to the properties of the maps in $\swi$. We will show next that 
$v$ can be extended to $[0,1)$ in such a way that $v \in \swi$. 
(In the case of the map defined by (6) we take
\begin{equation}\label{eq:f}
\dst v(x)=\left \{       
\begin{array}{l}               
\dst \frac{x+3}{4} \quad \text{on} \quad 
\left[0,\frac{1}{4}\right),\\ \\
\dst \frac{x+1}{2}  \quad \text{on} \quad 
\left[\frac{1}{4}, \frac{1}{2}\right)\cup \left[\frac{5}{8}, \frac{3}{4}\right),\\ \\
\dst x \quad \text{on} \quad 
\left[\frac{1}{2},\frac{5}{8}\right)\cup \left[\frac{7}{8},1\right) ,\\ \\
\dst \frac{x}{2}  \quad \text{on} \quad 
\left[\frac{3}{4}, \frac{7}{8}\right),
\end{array} \right.                                          
\end{equation}  
which is just one of the various possible extentions.)\par

Getting back to the general situation, it turns out that independently of what extensions one might consider for $u$ and $v$, the map $u\diamond v$ constructed as in (\ref{cabe}) is the same as $h$ a.e. The following lemma solves the existence of the extension $v$.
\begin{lemma}\label{lem:exv}
The map $v:\ca R_2\ra [0,1)$ defined by $v(x)=h^{-1}(x)$ $(x\in \ca R_2)$ can be extended to a map in $\swi$.
\end{lemma}  
\proof. \  Let us define the extension inductively in the following way. First we just pick a bijection $\phi:\NN\ra \{(k,l):k,l\in \NN\cup \{0\}\}$. For $n\in \NN$ we denote by $\phi_1(n)$ [resp. $\phi_2(n)$] the first component [resp. second component] of $\phi(n)$.
Then the initial step is to extend $v$ to $E_1:=\ca R_2\cup F_1$ where 
$F_1:=\{x\in [0,1)\setminus \ca R_2:
(x+\phi_1(1))/2^{\phi_2(1)}\in [0,1)\setminus v(\ca R_2)\}$. (Clearly $F_1$ may be empty.) In any case we define $v_1(x)=(x+\phi_1(1))/2^{\phi_2(1)}$ for all $x\in F_1$ and $v_1(x)=v(x)$ for $x\in \ca R_2$. Suppose we defined $v_n$ on $E_n:=E_{n-1}\cup F_n$ as an extension of $v_{n-1}$. Let $v_{n+1}$ be the extension of $v_n$ to $E_{n+1}:=E_n\cup F_{n+1}$ where $F_{n+1}:=\{x\in [0,1)\setminus E_n:
(x+\phi_1(n+1))/2^{\phi_2(n+1)}\in [0,1)\setminus v_n(E_n)\}$ defined as 
$v_{n+1}(x)=(x+\phi_1(n+1))/2^{\phi_2(n+1)}$ for all $x\in F_{n+1}$.
By the way these extensions are constructed it is easy to see that 
each $v_n$ is a one-to-one map.
If for some $n\in \NN$ we have $E_n=[0,1)$ a.e., then the proof is finished since $v_n$ would be the extension we were looking for. \par
We may assume then that $\mu([0,1)\setminus E_n)>0$ for every $n\in \NN$. In this case we
let $E_\infty=\cup_{n=1}^{\infty}E_n$ and define the extension $\wt v$ of $v$ to $E_\infty$ in the usual way: 
$\wt v(x)=v_n(x)$ if $x\in E_n$. Because $v_{n+1}$ is an extension of $v_n$ the map 
$\wt v$ is well defined. We claim that $E_\infty=[0,1)$ a.e. in which case then 
$\wt v\in \swi$ is the extension we want.\par
In order to prove this claim we proceed by way of contradiction and assume that 
$\mu(U)>0$ where $U:=[0,1)\setminus E_\infty$. One can easily show, using similar arguments to those in the proof of Lemma~\ref{lem:firstm}, that $\mu(\wt v(E_\infty))\le \mu(E_\infty)$. This shows that $\mu(V)>0$ where $V:=[0,1)\setminus \wt v(E_\infty)$.
Obviously, $U$ and $V$ are measurable sets. It is known (see \cite{re}) that 
the transformation $T(x)=\lfloor 2x \rfloor$ is an ergodic 
transformation on $[0,1)$ with respect to an invariant measure
which is equivalent to Lebesgue measure. This implies that $\mu(T^k(V)\cap U)>0$ for some 
$k\in \NN$. Equivalently we have $\mu((2^kV-l)\cap U)>0$ for some $l\in \NN\cup\{0\}$.
Let us write $\wt U$ for the set  $(2^kV-l)\cap U$. We have $(x+l)/2^k\in V$ for every $x\in \wt U$. 
If we let $n:=\phi^{-1}(k,l)$ it follows that $\wt U\subset E_n\cap U$ which contradicts the fact that $E_n\cap U=\emptyset $ ($E_n\subset E_\infty= [0,1)\setminus U$).\eproof

Returning to the proof of the Theorem~\ref{main} let $u\in \fwi$ and $v\in \swi$ be the maps constructed as above. We need to show that $u\diamond v=h$ almost everywhere. Since 
$\ca R_1:=u(\ca D_1)$ and $u_{|\ca D_1}=h_{|\ca D_1}$ it follows that 
$h^{-1}(\ca R_1)=\ca D_1$ and $h^{-1}(\ca R_2)=\ca D_2$ since $h$ is one-to-one. Hence, $v({\ca R_2})=h^{-1}(\ca R_2)=\ca D_2$. Since $u$ and $v$ are one-to-one maps, it follows that $u({\ca D_2})\subset \ca R_2$ and $v({\ca R_1})\subset \ca D_1$.
Let $\ca S=[0,1)\setminus v([0,1))$. Clearly $\ca S\subset \ca D_1$. Thus, $u({\ca S})\subset \ca R_1$ and then $(v\circ u)({\ca S})\subset \ca D_1$. Inductively it follows that all the sets $(v\circ u)^k({\ca S})$, $k\in \Bbb N\cup {0}$, are contained in $\ca D_1$. Therefore for $x\in \bigcup_{k\ge 0} (v\circ u)^k({\ca S})$ according to (\ref{cabe}) we have $u\circ v(x)=u(x)=h(x)$.\par
Let $\ca N=v([0,1))\setminus (v\circ u)([0,1))$. Clearly $\ca N\subset \ca R_2$ and inductively it follows that all the sets $(v\circ u)^k({\ca N})\subset \ca R_2$. Thus,
for $x\in \bigcup_{k\ge 0} (v\circ u)^k({\ca N})$ we have 
$u\diamond v(x)=v^{-1}(x)=(h^{-1})^{-1}(x)=h(x)$. \par
 Using Lemma~\ref{lem:firstm} we observe that the set $\bigcap _{k\ge 0} (v\circ u)^k([0,1))$ must have Lebesgue measure zero. Hence $u\diamond v=h$ almost everywhere.\eproof
\n {\bf Remark.} If $v=id$ then $u\diamond v=id$. We will use this fact to obtain as a corollary Speegle's result in \cite{sp}. The convergence seems to be in a stronger metric but in fact the two metrics are equivalent on the class of wavelets.\par
\begin{theorem}\label{connect} For every two wavelet sets $W_0$ and $W_1$ one can find a chain of wavelet sets $\{W_t\}_{t\in [0,1]}$ connecting them and such that $d(W_t,W_s)\rightarrow 0$ if $t\rightarrow s$ where the metric $d$ is given by $(\ref{distance})$. 
\end{theorem}
\proof. It is obvious that we just need to consider the case $W_1=E$ (the Littlewood-Paley wavelet set). Denote $W_0$ simply by $W$ and let $h:=h_W$ be the wavelet induced isomorphism of $W$. 
According to the Theorem~\ref{main} we can find $u\in \fwi$ and $v\in \swi$
such that $u\diamond v=h$ almost everywhere. The idea of our proof is to connect $v$ with $id$ by a continuous chain  of maps in $\swi$ and then use the second part of Theorem~\ref{main} to construct 
$h_t=u\diamond v_t$. Then we just take $W_t$ the corresponding wavelet set to $h_t$.
\begin{lemma} \label{vt}
There exist a chain $\{v_t\}_{t\in [0,1]}$ such that $v_0=v$, $v_1=id$ and 
$\mu(\{x:v_t(x)\not=v_s(x)\})\rightarrow 0$ as $t\rightarrow s$. 
\end{lemma} 
\proof. We may assume without loss of generality that $v$ is not the identity function. Hence there must exist a measurable set $U\subset [0,1)$ such that $\mu(v(U))< \mu(U)$. Since $\mu(v(L))\le \mu(L)$ for every measurable subset of $[0,1)$ we see that the set $V=[0,1)\setminus v([0,1))$ has positive Lebesgue measure. We define what is going to be the first part of our chain by $v^1_t(x)=v(x)$ if $x<t$ and $x\in [0,1)\setminus V$ or $x\ge t$ and 
$v^1_t(x)=x$ if $x<t$ and $x\in V$. Clearly if $0\le t<s\le 1$ then $\{x:v^1_t(x)\not=v^1_s(x)\}\subset [t,s]$ and therefore 
$\mu(\{x:v^1_t(x)\not=v^1_s(x)\}\rightarrow 0$ as $t\rightarrow s$. 
Notice that the maps $v^1_t$ are by construction one-to-one 
and then automatically $v^1_t\in \swi$ for all $t$. The map $v^1_1$ is
$id$ on $V$ and $v$ everywhere else. Notice that $[0,1)\setminus(v_1^1([0,1))=v(V)$. Let us denote by $\ca Z$ the set of fixed points of $v$.
We have clearly $V\cap \ca Z=\emptyset$ and therefore $v^k(V)\cap \ca Z =\emptyset$ for all $k\in \Bbb N$. We observe that $[0,1)$ can be partitioned as $[0,1)=\bigcup_{k\ge 0} v^k(V)\cup \bigcap_{k\ge 0}v^k([0,1))$. Using the above observations 
$\ca Z \subset \bigcap_{k\ge 0}v^k([0,1))$ and since $\displaystyle \mu(v(L))\le \frac{\mu(L)}{2}$
for every measurable set contained in $[0,1)\setminus \ca Z$ we conclude that
$\bigcap_{k\ge 0}v^k([0,1))\setminus \ca Z$ has Lebesgue measure zero. Hence to finish the 
proof we will extend the chain $v^1_t$  with a chain $v^2_t$ which will connect $v^1_1$ to the $v^2_1$ which is $id$ on $V\cup v(V)$. 
So, the next chain is defined by $v^2_t(x)=v^1_1(x)$ if $x<t$ and $x\in [0,1)\setminus v(V)$ or $x\ge t$ and $v^2_t(x)=x$ if $x<t$ and $x\in v(V)$. Observe that maps $v^2_t$ are one-to-one and so they are in $\swi$. As before the continuity of $\{v^2_t\}_t$ is insured. The map $v^2_1$ is $id$ on $V\cup v(V)$ and $v$ everywhere else. Then we 
continue inductively constructing $v^3_t$, $v^4_t$,..., in a similar manner. To end the proof we put these countably many chains together in an obvious way to form the required chain. More precisely, we scale and glue the chains $\{v_t^1\}, \{v_t^2\}, \{v_t^3\},...$ to create the final chain $\{v_t\}$
where $v_1=id$ (we allocate an interval of length $1/2$ for $\{v_t^1\}$, an interval of length $1/4$ for $\{v_t^2\}$ and so on). It is easy to see that the continuity claim still holds, i.e., that 
$\mu(\{x:v_t(x)\not=v_s(x)\})\rightarrow 0$ as $t\rightarrow s$. \eproof

Returning to the proof of Theorem~\ref{connect} we want to establish next the equivalent of the distance given by (\ref{distance}) at the level of $\wi$.
\begin{lemma}\label{lem:distancewi} Let $h_1$ and $h_2$ be the wavelet induced isomorphisms associated with two wavelet sets $W_1$ and $W_2$ then the metric given by $(\ref{distance})$ satisfies:
\begin{equation}\label{distancewi}
d(W_1,W_2)=(4\pi\mu(\omega'))^{1/2}+(2\nu(\omega)\}))^{1/2},
\end{equation}
where the measure $\nu$ on $[0,1)$ is given by $\ds d\nu(x)=\frac{1}{1-x}\chi_{[0,1/2)}(x)d\mu(x)+\frac{1}{x}\chi_{[1/2,1)}(x)d\mu(x)$ and 
$\omega':=\{x\in [0,1):h^{-1}_1(x)\not=h^{-1}_2(x)\}$ and 
$\omega:=\{x\in [0,1):h_1(x)\not=h_2(x)\}$.
\end{lemma}
\proof.\  In order to establish (\ref{distancewi}) we observe that the property (iv) in Proposition~\ref{prop:charh}
implies that $W_1\bigtriangledown W_2=
\left(\delta_{|W_1}^{-1}\circ \xi^{-1}([0,1))\setminus \delta_{|W_1}^{-1}\circ \xi^{-1}(\Omega)\right) \bigcup 
\left(\delta_{|W_2}^{-1}\circ \xi^{-1}([0,1))\setminus \delta_{|W_1}^{-1}\circ \xi^{-1}(\Omega)\right)$ (disjoint union). Hence 
\begin{equation}\nonumber
\begin{array}{c}
\mu(W_1\bigtriangledown W_2)=\mu\left(\delta_{|W_1}^{-1}\circ \xi^{-1}([0,1)\setminus \Omega)\right) +
\mu\left(\delta_{|W_2}^{-1}\circ \xi^{-1}([0,1))\setminus \Omega \right)=\\
\mu\left(\delta_{|W_1}^{-1}\circ \xi^{-1}(\omega)\right)+
\mu\left(\delta_{|W_2}^{-1}\circ \xi^{-1}(\omega )\right).
\end{array}
\end{equation}
Taking in account that $\mu$ is invariant under translations and homogeneous under dilations, i.e., $\mu(tU)=t\mu(U)$ for every measurable set $U$ and every positive real number $t$, we obtain  
\begin{equation}\label{muw1w2}
\begin{array}{c}
\mu(W_1\bigtriangledown W_2)=\mu\left(\tau_{|W_1}\circ \delta_{|W_1}^{-1}\circ \xi^{-1}(\omega)\right)+
\mu\left(\tau_{|W_2}\circ \delta_{|W_2}^{-1}\circ \xi^{-1}(\omega)\right)=\\
2\pi\mu\left(\xi\circ\tau_{|W_1}\circ \delta_{|W_1}^{-1}\circ \xi^{-1}(\omega)\right)+
2\pi\mu\left(\xi\circ\tau_{|W_2}\circ \delta_{|W_2}^{-1}\circ \xi^{-1}(\omega)\right)=\\
2\pi\mu\left(h_1(\omega)\right)+
2\pi\mu\left(h_2(\omega)\right).
\end{array}
\end{equation}
Let us observe that $h_1(\omega)=h_2(\omega)=\omega'$ and so $h^{-1}(\omega')=h^{-1}(\omega')=\omega$. Indeed, $h_1(\omega)=h_1([0,1)\setminus \Omega)=[0,1)\setminus h_1(\Omega)$ and one can easily check that $h_1(\Omega)=\Omega'$ where $\Omega'=\{x:h_1^{-1}(x)=h_2^{-1}(x)\}$ where for the convenience of the reader we recall that 
$\Omega=\{x\in [0,1):  h_1(x)=h_2(x)\}$. As a result, 
(\ref{muw1w2}) becomes
\begin{equation}\label{dist1}
\mu(W_1\bigtriangledown W_2)=4\pi\mu\left(\omega'\right).
\end{equation}
This equality allows us to get the first term in (\ref{distancewi}). In order to
obtain the second part of (\ref{distancewi}) we begin by considering the measure on 
$\RR \setminus \{0\}$ 
defined as $d\lambda(x)=\frac{1}{|x|}d\mu(x)$. As before we have 
\begin{equation}\nonumber
\begin{array}{c}
\lambda(W_1\bigtriangledown W_2)=
\lambda\left(\delta_{|W_1}^{-1}\circ \xi^{-1}(\omega)\right)+
\lambda\left(\delta_{|W_2}^{-1}\circ \xi^{-1}(\omega) \right).
\end{array}
\end{equation}
But $\lambda$ is invariant under dilations.  Hence the above changes into:
\begin{equation}\label{dist2}
\begin{array}{c}
\lambda(W_1\bigtriangledown W_2)=
2\lambda\left(\xi^{-1}(\omega)\right).
\end{array}
\end{equation}
One can easily check that $\lambda\circ\xi^{-1}=\nu$. Putting (\ref{dist1}), (\ref{dist2}), and (\ref{distance}) together we obtain (\ref{distancewi}).\eproof
Returning now to the proof of Theorem~\ref{connect} let us define $\{h_t\}_{t\in [0,1]}$ by  $h_t:=u\diamond v_t$, with $v_t$ given by Lemma~\ref{vt}. According to Theorem~\ref{main} and Proposition~\ref{prop:charh} each $h_t$ is an wavelet induced isomorphism which corresponds to some  wavelet $W_t$. Because $h_0=u\diamond v$ and $h_1=u\diamond id$ we have $W_0=W$ and $W_1=E$. We need to introduce one more
notation: for two maps $f$ and $g$ having the same domain of definition let $\omega(f,g)$ be defined by
\begin{equation}\label{eq:omegafg}
\omega(f,g)=\{x:f(x)\not=g(x)\}).
\end{equation} 
By the way $\{v_t\}$ was constructed if $0\le t<s\le 1$ then $\omega(v_t,v_s)$ is contained in an interval of length less than $s-t$. By the definition (\ref{cabe}) we see that  $\omega(h_t,h_s)$ is contained in a set of Lebesgue measure less than $2[(t-s)/2+(t-s)/4+...+(t-s)/2^k+...]=2(t-s)$ because $u$ is ``measure contractive" is the sense of Lemma~\ref{lem:firstm}. Therefore $\mu(\omega(h_t,h_s))\rightarrow 0$ as $t\rightarrow s$ or $s\rightarrow t$. To finish the proof we observe that $\nu$ is equivalent with the Lebesgue measure on $[0,1)$ and $h_t^{-1}=v_t\circ u$. Using the same arguments as above we have $\mu(\omega( h^{-1}_t, h^{-1}_s))\rightarrow 0$ as 
$t\rightarrow s$ or $s\rightarrow t$. Finally we use Lemma~\ref{lem:distancewi} to end the proof.\eproof
{\bf Remark.} The Theorem~\ref{connect} has the advantage of being a more constructive result than the one in
\cite{sp}. We observe also that by construction $W_t\subset W\cup E$ for every $t$. This construction is in some sense very similar to the one given in \cite{ilp} where only a partial result was obtained.

\begin{theorem}[Garrigos-Speegle\cite{gs}]\label{compl} The class $\ws$ is complete in the metric
given by $(\ref{distance})$.
\end{theorem}

\proof. Let us start with a sequence of wavelet set $\{W_n\}$ which is Cauchy in the metric
in $(\ref{distance})$. Let us consider their corresponding wavelet induced isomorphisms
$\{h_n\}_{n\in \Bbb N}$. According to Lemma~\ref{lem:distancewi}, we have $\mu(\omega(h_n,h_m))\rightarrow 0$ and $\mu(\omega(h^{-1}_n,h^{-1}_m))\rightarrow 0$
as $m,n\rightarrow \infty$. Since $\int_0^1 |h_n(x)-h_m(x)|^2d\mu(x)\le \mu(\omega(h_n,h_m))$ 
and $\int_0^1 |h_n^{-1}(x)-h_m^{-1}(x)|^2d\mu(x)\le \mu(\omega(h_n^{-1},h_m^{-1}))$ it follows that $h_n$ and $h_n^{-1}$ are Cauchy sequences in $L^2([0,1))$. Let $f$, $g$  be their limits in $L^2([0,1))$. Passing to a subsequence we may assume that these sequences are pointwise convergent a.e.($[\mu]$). 
For each $n$ let $[0,1):=U_n\cup V_n$ be the decomposition (partition) which corresponds to the way we defined the maps $u_n$ and $v_n$ as in the Theorem~\ref{main} for the wavelet induced map $h_n$. In other words, $u_n(x)=h_n(x)$ for $x\in U_n$ and $v_n(y)=h_n^{-1}(y)$ for $y\in h(V_n)$. We remind the reader that $u_n$ and $v_n$ can be extended in order that $u_n\in \fwi$ and $v_n\in \swi$. \par
For $m>n$ we observe that $(U_n\bigtriangledown U_m)\setminus \ca F \subset \omega(h_n,h_m)$ where $\ca F$ was defined
by (\ref{calf}) in the proof of Proposition~\ref{prop:charh}. Indeed, for instance if $x\in (U_n\setminus U_m)\setminus \ca F =U_n\cap V_m$ it follows that $h_n(x)=u_n(x)$ and $h_m(x)=v_m^{-1}(x)$. Since $x\not\in \ca F$ we see that $h_n(x)\not=h_m(x)$. Hence $\mu(U_n\bigtriangledown U_m)\rightarrow 0$ as $m,n \rightarrow \infty$. Passing again to a subsequence if necessary we may assume without loss of generality that $\chi_{U_n}\rightarrow \chi_U$ pointwise for some measurable subset $U$ of $[0,1)$. It follows then that $\mu(U_n\bigtriangledown U)\rightarrow 0$,  $\mu(V_n\bigtriangledown V)\rightarrow 0$ as $n\ra \infty$ and $\chi_{V_n}\rightarrow \chi_V$ pointwise where $V=[0,1)\setminus U$. It is clear that $U$ may be empty in which case we will show that $f=g=id$. If $\mu(U)>0$ then $\chi_{U_n}(x)\ra \chi_U(x)=1$ for almost every $x\in U$. For such an $x\in U$ we deduce that for some $n_0(x)\in \Bbb N$
if $n\ge n_0(x)$ it follows that $x\in U_n$. Hence $h_n(x)=u_n(x)=\lfloor x/2^{k(x,n)}\rfloor\ra f(x)$ or $h_n(x)=u_n(x)=\lfloor (x-1)/2^{k(x,n)}\rfloor \ra f(x)$ with $k(x,n)\in \Bbb N$. There are only two possibilities: $k(n,x)\ra \infty$ or $k(n,x)$ is eventually constant. We claim that 
for almost every $x\in U$ the sequence $k(x,n)$ is eventually constant. \par
Indeed, suppose that for some subset of $U$ of positive measure, say $\Upsilon$, 
we have $k(x,n)\ra \infty$ for all $x\in \Upsilon$. Then if we fix $n$ but keep it arbitrary  it is clear from our assumption on $\Upsilon$ that $\Upsilon=\bigcup_m\{x\in \Upsilon:\ k(x,m)>k(x,n)\}$. Using the continuity of the Lebesgue measure we may replace $\Upsilon$ with  $\Upsilon_{m,n}=\{x\in \Upsilon:\ k(x,m)>k(x,n)\}$ for some $m$ such that $\mu(\Upsilon_{m,n})>\mu(\Upsilon)/2$. Then it easy to see that $\Upsilon_{m,n}\subset \omega(h_n,h_m)$ which implies $0<\mu(\Upsilon)/2<\mu(\omega(h_n,h_m))$. This contradicts the fact $\mu(\omega(h_n,h_m))\ra 0$ and our claimed is proved.\par
We proved that for almost every $x\in U$ it follows that $k(x,n)$ is an eventually constant
sequence. For such an $x$ we let $\ds k(x)=\lim_{n\ra \infty}k(x,n)\in \Bbb N$. Then 
$f(x)=\lim h_n(x)=u_n(x)=\lfloor x/2^{k(x)}\rfloor$ or $f(x)=\lfloor (x-1)/2^{k(x)}\rfloor$ depending upon $x\in [1/2,1)$ or $x\in [0,1/2)$. By Lemma~\ref{lem:firstm}, $f$ is automatically one-to-one on $U$ and agrees with the characterization of an wavelet induced map in Proposition~\ref{prop:charh}.\par
We consider $U_n'=h_n(U_n)$ and $V_n'=h_n(V_n)$. Since $[0,1)=U_n\cup V_n$ we get $[0,1)=U_n'\cup V_n'$. As before for $m>n$ we observe that $(V_n'\bigtriangledown V_m')\setminus \ca F \subset \omega(h_n^{-1},h_m^{-1})$. Then without loss of generality
we may assume that $\chi_{V_n'}\rightarrow \chi_{V'}$ pointwise for some measurable subset $V'$ of $[0,1)$. As a result $\chi_{U_n'}\rightarrow \chi_{U'}$ pointwise a.e. where 
$U'=[0,1)\setminus V'$. For a.e. $x\in U$ we proved that there exist an $n_1(x)\in \Bbb N$ such that if $n>n_1(x)$ we have $f(x)=h_n(x)$ and $x\in U_n$. Since $h_n(x)\in h_n(U_n)=U_n'$ we get $\chi_{U_n'}(f(x))=1$. Then letting $n\ra \infty$ we get $\chi_{U'}(f(x))=1$. This shows that $f(U)=U'$ almost everywhere. Therefore if $\mu(V')=0$ then it must be true that $f(U)=[0,1)$ a.e. ($[\mu]$). But this is not possible because $f$ is strictly contractive in the sense of measure.\par
Thus, it must be true that $\mu(V')>0$. As before let us take an $y\in V'$ such that 
$\chi_{V_n'}(y)\rightarrow \chi_{V'}(y)=1$. Thus there exists an $n_2(y)\in \Bbb N$ such that $y\in V_n'$ if $n>n_2(y)$. By Proposition~\ref{prop:charh} and the construction of $u_v$, $v_n$ in Theorem~\ref{main} we have $h_n^{-1}(y)=v_n(y)=(y+k(y,n))/2^{l(y,n)}$ for some integer $k(y,n)$ and $l(y,n)\in \Bbb N\cup \{0\}$. We now claim that for almost every $y\in V'$ the sequence $l(y,n)$ is eventually constant on a subsequence. \par
Indeed, assuming that for some measurable subset of $V'$, say $\Theta$, with $\mu(\Theta)>0$, for $y\in \Theta$ the sequence $l(y,n)\rightarrow \infty$. Then we make the argument as above that for fixed $n$ but arbitrary  $\Theta=\bigcup_m\{y\in \Theta:\ l(y,m)>l(y,n)\}$. Hence, there exists $\Theta_{m,n}:=\{x\in \Upsilon:\ k(x,m)>k(x,n)\}$ for some $m$ such that $\mu(\Theta_{m,n})>\mu(\Theta)/2$. One can check similarly that $\Theta_{m,n}\setminus \ca F\subset \omega(h_n^{-1},h_m^{-1})$ which implies $0<\mu(\Theta)/2<\mu(\omega(h_n^{-1},h_m^{-1}))$. This contradicts the fact $\mu(\omega(h_n^{-1},h_m^{-1}))\ra 0$ and our claimed is proved. \par
Using this claim and the fact that $h_n^{-1}(y)$ is convergent for almost every $y$ to $g(y)$ we see that the sequence $k(y,n)$ must be eventually constant for a subsequence
on which $l(y,n)$ is eventually constant. Letting $n\ra \infty$ we obtain that $g(y)=(y+k(y))/2^{l(y)}$ for some 
$k(y),\ l(y)\in \Bbb N\cup \{0\}$. This shows that $g$ is 
in $\swi$. As we argued before $g(V')=V$. If $\mu(U)=0$ then $\mu(U')=0$ and then $g$ must be the identity since it is contractive in the sense of measure. \par
To finish the prove we extend $f$ from $U$ to $[0,1)$ and $g$ from $V$ to $[0,1)$ so that 
the new maps $\wt f$ and $\wt g$ satisfy: $\wt f\in \fwi$ and $\wt g\in \swi$. Denote $h=\wt f\diamond \wt g$. Finally an argument based on the continuity of the Lebesgue measure and the relation (\ref{cabe}) shows that $h_n$ converges to $h$.\eproof  
{\bf Remark.} The proof of Theorem~\ref{compl} although quite lengthy it reveals the understanding of what the convergence in the distance (\ref{distance}) means at the level of wavelet sets and equivalently on $\wi$.

\par
\section{Examples}
Consider the Journe wavelet set $\dst J:=\left[-\frac{32\pi}{7},-4\pi\right)\cup \left [-\pi,-\frac{4\pi}{7}\right)\cup \left
[\frac{4\pi}{7},\pi\right)\cup \left [4\pi,\frac{32\pi}{7}\right)$. Then its wavelet induced function denoted by $\wt h_J$ can be described by
\begin{equation}\label{eq:journe}
\dst \wt h_J(x)=\left \{       
\begin{array}{l}               
\dst \lfloor 2^{-1}(x-1) \rfloor \quad \text{on} \quad 
\left[0,\frac{3}{7}\right),\\ \\
\dst \lfloor 2^2x \rfloor  \quad \text{on} \quad 
\left[\frac{3}{7}, \frac{4}{7}\right),\\ \\
\dst \lfloor 2^{-1}x \rfloor  \quad \text{on}  \quad              
\left[\frac{4}{7}, 1\right).
\end{array} \right.                                          
\end{equation}  
In this case one can choose $u\in \fwi$ defined by
\begin{equation}\label{eq:jpurneu}
\dst u(x)=\left \{       
\begin{array}{l}               
\dst \lfloor 2^{-1}(x-1) \rfloor \quad \text{on} \quad 
\left[0,\frac{1}{2}\right),\\ \\
\dst \lfloor 2^{-1}x \rfloor  \quad \text{on}  \quad              
\left[\frac{1}{2}, 1\right),
\end{array} \right.                                          
\end{equation}  
and $v\in \swi$ defined by
\begin{equation}\label{eq:jpurnev}
\dst v(x)=\left \{       
\begin{array}{l}               
\dst \frac{x+2}{4}  \quad \text{on} \quad 
\left[0,\frac{1}{2}\right),\\ \\
\dst \frac{x+1}{4}   \quad \text{on}  \quad              
\left[\frac{1}{2}, 1\right),
\end{array} \right.                                          
\end{equation}  
\bigskip 
such that $\wt h_J=u\diamond v$.\par

Let us go the other way in our construction. 
Arguably, the simplest map $u\in \fwi$ that one can take is given in (\ref{eq:jpurneu}) and we pick one of the simplest maps $v\in \swi$ defined by $v(x)=x/2$ for all $x\in [0,1)$. Then the wavelet induced map, $h=u\diamond v$, that is obtained by the Schr\"{o}der-Cantor-Bernstein construction is defined by
\begin{equation}\label{eqex1}
\dst h(x)=\left \{       
\begin{array}{l}               
\dst \frac{x}{2}  \quad \text{on} \quad 
\left[\frac{1}{2},1\right),\\ \\
\dst 2x \quad \text{on} \quad 
\left[\frac{1}{3},\frac{1}{2}\right),\\ \\
\dst \frac{x+1}{2} \quad \text{on} \quad 
A,\\ \\
\dst 2x \quad \text{on} \quad 
\left[0,\frac{1}{3}\right)\setminus A,
\end{array} \right.                                          
\end{equation}  
where $A=\dst \bigcup_{n=1}^{\infty}[z_n,x_{n+1})$ and the sequences $\{z_n\}_{n\ge 1}$, $\{x_n\}_{n\ge 1}$ are given by the formulae 
$\dst z_n=\frac{1}{3}\left(1-\frac{5}{2\times 4^{n}}\right)$, $\dst x_n=\frac{1}{3}\left(1-\frac{1}{4^{n-1}}\right)$. This will give rise to a wavelet set containing infinitely many intervals $\dst W:=\left[\frac{\pi}{2},\pi\right)\cup 
\left[\frac{-8\pi}{3},-2\pi\right)
\cup (\pi A-\pi)\cup \left([-4\pi,\frac{-8\pi}{3})
\setminus (4\pi A-4\pi)\right)$.
It is interesting to mention that if one modifies $v$, as in the proof of Theorem~\ref{connect}, to $v(x)=x/2$ for $x\in [0,2/7)$ and $v(x)=x$ if $x\in [2/7,1)$
then the wavelet set constructed from $u$ and the new $v$ has only six intervals.

{\bf Acknowledgments.} First I would like to express my thanks to Professor David Larson for introducing me to the subject of wavelet theory years ago, and for his encouragements that determined me to pursue this research. Proposition~\ref{prop:charh} appears just with statement (parts (i)-(iii)) in \cite{ailp} and the seed of Theorem~\ref{main} was the subject of discussions with Professor E. Azoff while I was a postdoctoral associate at University of Georgia, Athens. Hence I am greatly thankful to Professor Azoff for his collaboration and continuous support during a period of more than five years.
Last but not least I am thankful to the refree who had the patience to carefully read a previous version of this paper and
who pointed out to me how to significantly improve the presentation of the material.

\vfill \eject
     
\end{document}